\documentclass[11pt]{amsart}

\usepackage{amssymb}
\usepackage{comment}
\usepackage[all]{xy}
\usepackage{hyperref}

\newcommand{\R}{\mathbb R}
\newcommand{\N}{\mathbb N}
\newcommand{\Z}{\mathbb Z}
\newcommand{\ds}{\displaystyle}
\newcommand{\lie}{\mathcal{L}}
\renewcommand{\d}{\mathrm d}

\newcommand{\ra}{\rightarrow}
\newcommand{\p}{\partial}
\newcommand{\br}{\left[\, , \, \right]}
\newcommand{\brr}[1]{\left[#1\right]}
\newcommand{\X}{\ensuremath{\mathfrak{X}}}

\newcommand{\al}{\alpha}

\newcommand{\ran}{\rangle}
\newcommand{\lan}{\langle}

\DeclareMathOperator{\m}{mod} \DeclareMathOperator{\tr}{Tr}

\newtheorem{theorem}{Theorem}           %    [section]
\newtheorem{prop}[theorem]{Proposition}

\newtheorem{corol}[theorem]{Corollary}
\newtheorem{defi}[theorem]{Definition}
\newtheorem{ex}[theorem]{Example}

\newtheorem{rmk}[theorem]{Remark}

\begin{document}

\title{Jacobi-Nijenhuis algebroids and their modular classes}
\author{Raquel Caseiro}
\address{CMUC, Department of Mathematics, University of Coimbra}\email{raquel@mat.uc.pt}
\author{Joana M. Nunes da Costa}
\address{CMUC, Department of Mathematics, University of Coimbra}\email{jmcosta@mat.uc.pt}
\keywords{Jacobi algebroid, Jacobi-Nijenhuis algebroid,  modular
class}

\subjclass[2000]{17B62, 17B66, 53D10, 53D17}

\thanks{This work was
partially supported by POCI/MAT/58452/2004 and CMUC/FCT}

\begin{abstract}
Jacobi-Nijenhuis algebroids are defined as a natural generalization
of Poisson-Nijenhuis algebroids, in the case where there exists a
Nijenhuis operator on a Jacobi algebroid which is compatible with
it. We study modular classes of Jacobi and Jacobi-Nijenhuis
algebroids.
\end{abstract}

 \maketitle

\section{Introduction}

It is well known that the cotangent bundle $T^\ast M$ of any Poisson
manifold $M$ admits a Lie algebroid structure  and the pair
$(TM,T^\ast M)$ is a Lie bialgebroid over $M$. As a kind of
reciprocal result, any Lie bialgebroid $(A,A^\ast)$ induces a
Poisson structure on its base manifold. A special kind of Lie
bialgebroids are the triangular Lie bialgebroids. These are Lie
algebroids $(A,P)$ equipped with an $A$-bivector field $P$ such that
$[P,P]=0$. The $A$-bivector field $P$ induces a Lie algebroid
structure on the dual vector bundle $A^\ast$ so that the pair
$(A,A^\ast)$ is a Lie bialgebroid. Triangular bialgebroids are also
called Lie algebroids with a Poisson structure.

If one moves from the Poisson to the Jacobi  framework, these
statements are not true. In fact, if $M$ is a Jacobi manifold, its
cotangent bundle $T^\ast M$ is not, in general, a Lie algebroid. In
order to associate a Lie algebroid to a Jacobi manifold, one has to
consider the $1$-jet bundle $T^\ast M \times \R \to M$. However, if
we take the dual vector bundle $TM \times \R \to M$ endowed with its
natural Lie algebroid structure, the pair $(TM \times \R, T^\ast M
\times \R)$ is not a Lie bialgebroid. Motivated by this, Iglesias
and Marrero \cite{IglMarr} and Grabowski and Marmo \cite{GraMar1}
introduced the concepts of Jacobi algebroid, i.e. a Lie algebroid
$A$ with a $1$-cocycle $\phi_0$,
 and
of Jacobi bialgebroid, i.e. a pair $((A,\phi_0), (A^\ast, X_0))$ of
Jacobi algebroids in duality satisfying a compatibility condition.
Jacobi bialgebroids  admit Lie bialgebroids as particular cases and
are well adapted to the Jacobi context since every Jacobi manifold
$(M,\Lambda,E)$ has an associated Jacobi bialgebroid, $((TM \times
\R,(0,1)), (T^\ast M \times \R, (-E,0)))$. Imitating the Poisson
case, Iglesias and Marrero introduced in \cite{IglMarr} the notion
of triangular Jacobi bialgebroid, as follows. If $(A,\phi_0)$ is a
Jacobi algebroid and $P$ is a Jacobi bivector field, i.e. an
$A$-bivector field such that $\brr{P,P}^{\phi_0}=0$, then there
exists a Lie algebroid structure on $A^\ast$ with a $1$-cocycle such
that the pair of Jacobi algebroids in duality is a Jacobi
bialgebroid.
 As it happens in the Poisson case,
 the base manifold of any Jacobi bialgebroid
inherits a Jacobi structure.

On the other hand, Poisson-Nijenhuis structures on Lie algebroids,
i.e. Poisson-Nijenhuis algebroids, were introduced by Grabowski and
Urbanski in \cite{GraUrb}, as Lie algebroids equipped with a Poisson
structure and a Nijenhuis operator fulfilling some compatibility
conditions. In the first part of this paper, we extend this concept
to the Jacobi framework and we study \emph{Jacobi-Nijenhuis
algebroids}.

The other main goal of this paper is to study modular classes,
including modular classes of Jacobi-Nijenhuis algebroids. The
modular class of a Poisson manifold was defined by Weinstein in
\cite{Weinstein}, as an analogue in Poisson geometry of the modular
automorphism group of a von Neumann algebra. In \cite{EvLuWe},
Evens, Lu and Weinstein introduced the notion of modular class of a
Lie algebroid $A$ over $M$, using a representation of $A$ on the
line bundle $\X ^{\mathrm{top}}(A) \otimes
\Omega^{\mathrm{top}}(M)$. For the case of the cotangent Lie
algebroid $T^\ast M$ of a Poisson manifold $M$, they showed that its
modular class is twice the modular class of $M$, in the sense of
\cite{Weinstein}.
 Modular classes of
triangular Lie bialgebroids were studied in \cite{Kosmann2}, from
the point of view of generating operators for Batalin-Vilkovisky
algebras.

Regarding the Jacobi context, the first work on modular classes is
due to Vaisman, who introduced in \cite{Vaisman} the concept of
modular class of a Jacobi manifold. Then, in \cite{IglLopMarPad},
modular classes of triangular Jacobi bialgebroids were studied.

In the second part of this paper, we consider a Jacobi algebroid
$(A, \phi_0)$ and we define, using the $1$-cocycle $\phi_0$, a new
representation of $A$ on the line bundle $\X ^{\mathrm{top}}(A)
\otimes \Omega^{\mathrm{top}}(M)$, which leads to the definition of
\emph{modular class of a Jacobi algebroid}.

Modular classes of Poisson-Nijenhuis algebroids were defined in
\cite{raq}. Inspired in \cite{raq}, we define \emph{modular class of
a Jacobi-Nijenhuis algebroid}. We obtain a hierarchy of vector
fields on the Jacobi algebroid that covers a hierarchy of Jacobi
structures on the base.

The paper is divided into five sections. Section 2 is devoted to
Jacobi algebroids. We recall how to obtain a Lie algebroid structure
on $A \times \R$ over $M\times \R$ from a Jacobi algebroid
$(A,\phi_0)$  over $M$ \cite{IglMarr,GraMar1}. The notion of
compatibility of two Jacobi bivectors on a Jacobi algebroid is
introduced, and we prove that these Jacobi bivectors cover two
compatible Jacobi structures on the base manifold. In Section 3 we
define Jacobi-Nijenhuis algebroid and we show that a
Jacobi-Nijenhuis algebroid defines a hierarchy of compatible Jacobi
bivectors on the Jacobi algebroid and a hierarchy of compatible
Jacobi structures on the base manifold. Moreover, the dual vector
bundle also inherits a hierarchy of Jacobi algebroid structures that
provides the existence of a family of triangular Jacobi
bialgebroids. As a particular case of this construction, we recover
the notion of strong (or strict) Jacobi-Nijenhuis manifold
\cite{IglMarr2,Cos}. In Section 4, we introduce the notion of
modular class of a Jacobi algebroid and we discuss the relation
between modular class of a Jacobi algebroid $(A, \phi_0)$ over $M$
and modular class of the Lie algebroid $A \times \R$ over $M \times
\R$. Relations between modular forms of $A^\ast$ and $A^\ast \times
\R$, in the triangular case, as well as duality between modular
classes of $A$ and $A^\ast$ are also discussed. At this point we
relate our results with those obtained in \cite{IglLopMarPad}. In
Section 5, we give the definition of modular class of a
Jacobi-Nijenhuis algebroid and we prove a result which generalizes
the corresponding one of \cite{raq}: there exists a hierarchy of
$A$-vector fields that defines two hierarchies of vector fields, one
on $M\times \R$ and another on $M$. These hierarchies determine a
family of Jacobi structures on the manifold $M$.

\bigskip

\noindent {\bf Notation and conventions:} Let $(A,\rho, \br)$ be a
Lie algebroid over $M$. We denote by $\X^k(A)$ (resp. $\Omega^k(A)$)
the $C^\infty(M)$-module of $A$-$k$-vector fields (resp.
$A$-$k$-forms), by $\X(A)=\oplus_k\X^k(A)$ (resp.
$\Omega(A)=\oplus_k\Omega^k(A)$) the corresponding Gerstenhaber
algebra of $A$-multivector fields (resp.  $A$-forms) and by
$\X^{\mathrm{top}}(A)$ the top-degree sections of $A$. The De Rham
differential is denoted by $d$  while $\d$ stands for the Lie
algebroid differential.

Regarding the conventions of sign for the Schouten bracket and for
the interior product by a multivector field, we use the same
conventions of \cite{GraMar1,Kosmann2}, which are different from
those of \cite{IglMarr, IglLopMarPad}.

\section{Jacobi algebroids}

We begin by recalling some well known facts about Jacobi algebroids.

\subsection{Jacobi algebroids} A \emph{Jacobi algebroid} \cite{GraMar1} or \emph{generalized
Lie algebroid} \cite{IglMarr} is a pair $(A, {\phi_0})$ where
 $A=(A, \br, \rho)$ is a Lie algebroid over a manifold $M$ and
 ${\phi_0}\in\Omega^1(A)$ is a 1-cocycle in the Lie algebroid cohomology with trivial coefficients, ${\ds \d {\phi_0}=0}$.
 A Jacobi algebroid has an associated
Schouten-Jacobi bracket on the graded algebra  $\X(A)$ of
multivector fields on $A$ given by
\begin{equation}\label{Jacobibracket}
    \brr{P,Q}^{\phi_0}=\brr{P,Q}+(p-1)P\wedge i_{\phi_0} Q- (-1)^{p-1}(q-1)
    i_{\phi_0} P\wedge Q,
\end{equation}
for $P\in \X^p(A), Q\in\X^q(A)$.

This bracket $\br^{\phi_0}$ satisfies the following properties (in
fact it is totaly defined by them), with $X,Y\in\X^1(A)$,
$P\in\X^p(A)$, $Q\in\X^q(A)$ and $f\in C^\infty(M)$:
\begin{equation}\label{JacobiBracket1}
    \brr{X,f}^{\phi_0}=\rho^{\phi_0}(X) f,
\end{equation}
\begin{equation}\label{JacobiBracket2}
 \brr{X,Y}^{\phi_0}=\brr{X,Y},
\end{equation}
\begin{equation}\label{JacobiBracket3}
    \brr{P,Q}^{\phi_0}=-(-1)^{(p-1)(q-1)}\brr{Q,P}^{\phi_0},
\end{equation}
\begin{equation}\label{JacobiBracket4}
    \brr{P, Q\wedge R}^{\phi_0}=\brr{P,Q}^{\phi_0}\wedge R +
    (-1)^{(p-1)q}Q\wedge \brr{P,R}^{\phi_0}- (-1)^{p-1}i_{\phi_0} P \wedge Q\wedge
    R,
\end{equation}
\begin{align}\label{JacobiBracket5}
    (-1)^{(p-1)(r-1)}\brr{P,\brr{Q,R}^{\phi_0}}^{\phi_0}&+(-1)^{(q-1)(p-1)}\brr{ Q,\brr{R,
    P}^{\phi_0}}^{\phi_0} \nonumber\\
    &+(-1)^{(r-1)(q-1)}\brr{R,\brr{P, Q}^{\phi_0}}^{\phi_0}=0.
\end{align}

In property (\ref{JacobiBracket1}), $\rho^{\phi_0}$ is the
representation of the Lie algebra  $\mathfrak{X}^1(A)$ on
$C^\infty(M)$ given by
$$\rho^{\phi_0}(X)f=\rho(X)f+f\lan {\phi_0},X\ran.$$ The cohomology
operator $\d^{\phi_0}$ associated with this representation is called
\emph{${\phi_0}$-differential} of $A$ and is given by
\begin{equation}\label{phidiff}
 \d^{\phi_0} \omega= \d\omega+{\phi_0}\wedge
 \omega,\quad\omega\in\Omega(A).
\end{equation}
With the ${\phi_0}$-differential we can  define a
\emph{${\phi_0}$-Lie derivative}:
\begin{equation}\label{phiLiederivative}
    \lie^{\phi_0}_X\omega=i_X\,\d^{\phi_0} \omega+(-1)^{p-1}\d^{\phi_0}\, i_X\omega, \quad
    X\in\X^p(A), \omega\in\Omega(A).
\end{equation}

In \cite{IglMarr, GraMar} we can find a construction which allow us
to obtain a Lie algebroid over $M\times \R$ from a Jacobi algebroid
over $M$. This construction is very  useful when we speak about
Jacobi algebroids, in fact it contains the essence of philosophy
adopted in the proofs in this paper, so we will explain it now.

 Consider the natural vector bundle $\hat A=A\times \R$ over $M\times
 \R$. The sections of $\hat A$ may be seen as time-dependent sections of
 $A$ and this space is generated as a $C^\infty(M\times \R)$-module by the
 space of sections of $A$, which are simply the time-independent
 sections of $\hat A$.

The anchor
\begin{equation}\label{anchorhatA}
    \hat\rho(X)=\rho(X)+\lan {\phi_0},X\ran \frac{\p}{\p t}, \quad
    X\in\X^1(A),
\end{equation}
and the bracket defined by $\br$ for  time independent multivectors
\begin{equation}\label{brackethatA}
    \brr{X,Y}_{\hat A}=\brr{X,Y}, \quad X,Y\in\X(A),
\end{equation}
define a Lie algebroid structure on $\hat A$ that we call the
\textit{induced Lie algebroid structure from $A$ by ${\phi_0}$}. If
$\hat \d$ is the differential in $\hat A$, from (\ref{anchorhatA})
we get
\begin{equation}\label{exactphi}
    {\phi_0}=\hat \d t,
\end{equation}
which means that the 1-cocycle ${\phi_0}$ can be seen as an exact
1-form on $\hat A$.

%%%%%%%% GAUGING ON X(A) %%%%%%%%

Considering the gauging in $\X(A)$ defined by
\begin{equation*}
 \tilde X=e^{-(p-1)t}X,\quad X\in\X^{p}(A),
\end{equation*}
we have the following relation between the Lie bracket in $\X(\hat
A)$ and the Jacobi bracket (\ref{Jacobibracket}):
\begin{equation}\label{relgauging1}
    \brr{\tilde X,\tilde Y}_ {\hat A}=\widetilde{\brr{X,Y}^{\phi_0}}.
\end{equation}

Now consider a \emph{Jacobi bivector} on $A$, i.e., a bivector
$P\in\X^2(A)$ such that
\begin{equation}\label{jacobibivector}
    \brr{P,P}^{\phi_0}=0.
\end{equation}
>From relation (\ref{relgauging1}) we deduce that $\tilde P=e^{-t}P$
is a Poisson bivector on $\hat A$ and, consequently, it defines a
Lie algebroid structure over $M\times \R$  on $\hat A^\ast$ given by
\begin{equation}\label{dualstructure:hatAdual}
    \brr{\alpha,\beta}_{\tilde P}=\widehat{\lie}_{\tilde
    P^\sharp\alpha}\beta-\widehat{\lie}_{\tilde P^\sharp\beta}\alpha-\hat\d
    \tilde{P}(\alpha,\beta),
\end{equation}
\begin{equation}
\hat\rho_\ast(\alpha)=\hat\rho\circ\tilde P^\sharp(\alpha)
\end{equation}
where $\alpha,\beta\in\X^1(\hat A^\ast)$ and $\widehat{\lie}$ is the
Lie derivative in $\hat A$. In particular, for $\alpha,\,
\beta\in\X^1(A^\ast)$, we have
\begin{equation}\label{relgauging2}
    \brr{e^t\alpha,e^t\beta}_{\tilde P}=e^t(\lie^{\phi_0}_{P^\sharp \alpha}\beta - \lie^{\phi_0}_{P^\sharp\beta}\alpha-\d^{\phi_0}
    P(\alpha,\beta)).
\end{equation}
The Lie bracket
\begin{equation}\label{bracketdualA}
    \brr{\alpha,\beta}_{P}={\lie}^{\phi_0}_{P^\sharp \alpha}\beta - {\lie}^{\phi_0}_{P^\sharp\beta}\alpha-\d^{\phi_0}
    P(\alpha,\beta),
\end{equation}
together with the anchor
\begin{equation}\label{anchordualA}
    \rho_\ast=\rho\circ P^\sharp,
\end{equation}
endows  $A^\ast$ with a Lie algebroid structure over $M$.

The section on $A$, $X_0=-P^\sharp({\phi_0})$ is a 1-cocycle of
$A^\ast$, and so $(A^\ast, X_0)$ is a Jacobi algebroid. The pair
$((A,{\phi_0}), (A^\ast, X_0))$ is a special kind of Jacobi
bialgebroid called \emph{triangular Jacobi bialgebroid} and we will
denote it by $(A,{\phi_0},P)$.

Recall that a \emph{Jacobi bialgebroid}  (see \cite{IglMarr},
\cite{GraMar1}) is a pair of Jacobi algebroids in duality,  ${\ds
((A, {\phi_0}), (A^\ast, X_0))}$, such that $\d_\ast^{X_0}$ is a
derivation of $(\X(A), \br^{\phi_0})$ or, equivalently,
$\d^{\phi_0}$ is a derivation of $(\X(A^\ast), \br^{X_0}_\ast)$.

\

The relation (\ref{relgauging2}) can be generalized to multisections
of $A^\ast$ if we consider the gauging
 in $\Omega(A)$:
\begin{equation}\label{relgauging4}
    \hat \omega=e^{pt}\omega, \quad \omega\in\Omega^p(A).
\end{equation}

\begin{prop} Let $\alpha,\beta$ be multisections of $A^\ast$. Then
\begin{equation}\label{relgauging3}
    \brr{\hat\alpha,\hat \beta}_{\tilde P}=\widehat{\brr{\alpha,\beta}}_{P}.
\end{equation}
\end{prop}

One should also notice that the structure of Lie algebroid on $\hat
A^\ast$ does not coincide with Lie algebroid structure induced from
$A^\ast$ by  the 1-cocycle $X_0$ (at least not in the same way it
was done with $A$ and ${\phi_0}$). In fact, the bracket of two  time
independent sections on $\hat A^\ast$, $\al, \beta\in\Omega^1(A)$,
is given by
\begin{equation}\label{bracketchapeuX0}
\brr{\alpha,\beta}_{\tilde P}=e^{-t}(\brr{\alpha,\beta}_{P}-\lan
\alpha, X_0\ran \beta+\lan \beta, X_0\ran \alpha)
\end{equation}
and the anchor of $\hat A^\ast$ is defined by
$$\hat\rho_\ast(\alpha)=e^{-t}\left(\rho_\ast(\al)+\lan\al,X_0\ran\frac{\p}{\p t}\right).$$

Any  Jacobi bialgebroid $((A,{\phi_0}),(A^\ast,X_0))$ gives to $M$ a
structure of Jacobi manifold, i.e., it equips $M$ with a bivector
field $P_M$ and a vector field $E_M$ satisfying
\begin{equation}\label{Jacobi:manifold}
\brr{P_M,P_M}=-2E_M\wedge P_M, \quad \brr{E_M,P_M}=0,
\end{equation}
or, equivalently, it defines a Jacobi bracket on $C^\infty(M)$ given
by:
$$
\{f,g\}_M=\lan \d^{\phi_0} f,\d_\ast^{X_0} g\ran.
$$

In particular, if $(A,{\phi_0},P)$ is a triangular Jacobi
bialgebroid then $(P_M,E_M)$ is defined by \footnote{We denote by
$\rho ^p$ the  morphism $\rho ^p : \X^p(A) \to \X^p(M)$, given by
$\rho^p
P(\alpha_1,\ldots,\alpha_p)=P(\rho^\ast\alpha_1,\ldots,\rho^\ast\alpha_p)$,
with $\alpha_1,\ldots,\alpha_p \in \Omega^1(M)$. Since $\rho$ is a
Lie algebroid morphism, we have that
$\rho^{p+q-1}\brr{P,Q}=\brr{\rho^p P,\rho^q Q}$, with $P \in
\X^p(A)$ and $Q \in \X^q(A)$.}
\begin{align}
    &P_M(df,dg)=\rho^2P(df,dg)=P(\rho^\ast df,\rho^\ast dg)=P(\d f,\d
    g),\label{jacobibivector:manifold:one}\\
    &E_M=\rho\circ P^\sharp({\phi_0}).\label{jacobibivector:manifold:two}
\end{align}

%%%%%%%%%%%%%%%%%%%%%%%%%% EXEMPLO: THE JACOBI BIALGEBROID OF A JACOBI MANIFOLD %%%%%%%%%%%%%%%%%%%%%%%%%%%%%%%%%%%

\subsection{The triangular Jacobi  bialgebroid of a Jacobi manifold}

Let $(M,\Lambda,E)$ be a Jacobi manifold, i.e, a manifold equipped
with a bivector $\Lambda$ and a vector field $E$ such that
\begin{equation}\label{jacobimanifold}
\brr{\Lambda,\Lambda}=-2E\wedge \Lambda, \quad \brr{E,\Lambda}=0.
\end{equation}
 The
vector bundle $T^\ast M\times \R$ is endowed with a Lie algebroid
structure over $M$ \cite{KerSou}. The Lie bracket and the anchor are
defined by
\begin{align}
    \brr{(\al,f),(\beta,g)}&_{(\Lambda,E)}=(\lie_{\Lambda^\sharp\al} \beta-\lie_{\Lambda^\sharp\beta}\al- d(\Lambda(\alpha,\beta))+f\lie_E\beta
    -g\lie_E\al\nonumber\\
    &-i_E(\alpha\wedge \beta),
    \Lambda(\beta, \al) +\Lambda(\alpha,d g)-\Lambda(\beta,d
    f)+fE(g)-gE(f))\label{bracket:TdualMXR}
\end{align}
and
\begin{equation*}
    \widetilde{(\Lambda,E)}^\sharp(\alpha,f)=\Lambda^\sharp(\al)+fE.
\end{equation*}

 In this Lie algebroid the
differential is given by
\begin{equation*}
\d_\ast(X,Y)=(\brr{\Lambda,X}+k E\wedge X + \Lambda \wedge Y,
-\brr{\Lambda,Y}-(k-1)E\wedge Y + \brr{E,X}),
\end{equation*}
for $(X,Y)\in\X^k(M)\oplus \X^{k-1}(M).$ The section $X_0=(-E,0)$ is
a 1-cocycle of $T^\ast M\times\R$ and the $X_0$-differential is
\begin{align*}
\d_{\ast}^{(-E,0)}(X,Y)=&(\brr{\Lambda,X}+(k-1) E\wedge X + \Lambda
\wedge Y, \\ &-\brr{\Lambda,Y}-(k-2)E\wedge Y + \brr{E,X}),
\end{align*}
for $(X,Y)\in\X^k(M)\otimes \X^{k-1}(M)$.

Now consider the canonical vector bundle $TM\times R$ over $M$ with
its structure of Lie algebroid given by the Lie bracket
$$
\brr{(X,f),(Y,g)}=(\brr{X,Y}, X(g)-Y(f))
$$
and the anchor
$$
\rho(X,f)=X.
$$

The differential $\d$ of this Lie algebroid is
$$
\d(\alpha,\beta)=(d \alpha, -d \beta),\quad \al,\beta\in\Omega(M).
$$
Obviously, ${\phi_0}=(0,1)$ is a 1-cocycle of $TM\times \R$. The
$\phi_0$-differential is given by
$$
\d^{(0,1)}(\alpha,\beta)=(d \al, \alpha-d \beta), \quad
\al,\beta\in\Omega(M).
$$

A Jacobi bivector on the Jacobi algebroid $(TM\times\R, (0,1))$ is a
section $(\Lambda,E)$ on $\X^2(M)\oplus \X^1(M)$ such that
\begin{equation}\label{jacobibivector:condition:manifold}
    \brr{(\Lambda, E), (\Lambda,E)}^{(0,1)}=0.
\end{equation}
Since (\ref{jacobibivector:condition:manifold}) is equivalent to
(\ref{jacobimanifold}), $(\Lambda,E)$ defines a Jacobi structure on
the manifold $M$. Moreover, $(\Lambda, E)^\sharp(0,1)=(E,0)$, where
${\ds (\Lambda,E)^\sharp: T^\ast M\times \R\ra TM\times \R }$ is the
vector bundle morphism defined by $(\Lambda,
E)^\sharp(\al,f)=(\Lambda^\sharp \al+fE,-i_E\alpha)$.

The Lie algebroid structure $(\br_{(\Lambda,E)},
\widetilde{(\Lambda, E)}^\sharp)$ in $T^\ast M\times \R$ coincides
with the Lie algebroid structure defined by the Jacobi bivector
$(\Lambda,E)$.  In fact one can check that
\begin{align*}
    \brr{(\al,f), (\beta,g)}_{(\Lambda,E)}=&\lie^{(0,1)}_{(\Lambda,E)^\sharp (\al,
    f)}  (\beta,g)  -\lie^{(0,1)}_{(\Lambda,E)^\sharp(\beta,
    g)}(\al,f)\\ &-\d^{(0,1)}\left((\Lambda,E)((\al,f),(\beta,g))\right)
\end{align*}
and
\begin{align*}
   \widetilde{ (\Lambda, E)}^\sharp=\rho\circ (\Lambda, E)^\sharp.
\end{align*}

 So we may conclude that
the pair $\left(\,(TM\times \R,(0,1)), (T^\ast M\times
\R,(-E,0))\,\right)$ is a triangular Jacobi bialgebroid
\cite{IglMarr}. Moreover, the Jacobi structure induced on the base
manifold coincides with the initial one.

%%%%%%%%%%%%%%%%%%%%%%%%%%%%%%%%% Compatible Jacobi bivectors %%%%%%%%%%%%%%%%%%%%%%%%%%%%%%%%%

\subsection{Compatible Jacobi bivectors}

With the construction presented in the section 2.1 the notion of
compatible Jacobi bivectors appears naturally.

\begin{defi}
Let $(A,\phi_0)$ be a Jacobi algebroid. Two Jacobi bivectors $P_1$
and $P_2$ on $A$ are said to be \textbf{compatible} if
\begin{equation}\label{jacobi:bivectors:compatible}
    \brr{P_1,P_2}^{\phi_0}=0.
\end{equation}
\end{defi}

Due to relation (\ref{relgauging1}), compatible Jacobi bivectors
$P_1$ and $P_2$ on $A$ are obviously associated with compatible
Poisson bivectors on $\hat A$, ${\tilde P_1}=e^{-t}P_1$ and
$\tilde P_2=e^{-t}P_2$:
$$
\brr{\tilde P_1,\tilde P_2}_{\hat A}=0.
$$
Moreover, they cover compatible Jacobi structures on the base
manifold $M$. Recall that two compatible Jacobi structures on a
manifold $M$ (see \cite{Cos2}) is a pair of Jacobi structures
$(\Lambda_1,E_1)$ and $(\Lambda_2, E_2)$ such that
$(\Lambda_1+\Lambda_2,E_1+E_2)$ is also a Jacobi structure, or,
equivalently, they satisfy the following two conditions:
\begin{align*}
    &\brr{\Lambda_1,\Lambda_2}=-E_1\wedge
    \Lambda_2-E_2\wedge\Lambda_1,\\
    & \brr{E_1,\Lambda_2}+\brr{E_2,\Lambda_1}=0.
\end{align*}

\begin{theorem}\label{compatible:jacobi:bivectors} Let
 $P_1$ and $P_2$ be compatible Jacobi bivectors on a Jacobi algebroid
$(A,\phi_0)$. These bivectors cover two compatible Jacobi structures
on the base manifold $M$.
\end{theorem}

\begin{proof} By definition of the Schouten-Jacobi bracket
$\br^{\phi_0}$, the compatibility condition ${\ds
\brr{P_1,P_2}^{\phi_0}=0}$ is equivalent to
\begin{equation}%\label{}
    \brr{P_1,P_2}=-P_{1}^\sharp(\phi_0)\wedge P_2-P_2^\sharp(\phi_0)\wedge P_1.
\end{equation}

On another hand, as we have mentioned, compatible Jacobi bivectors
$P_1$ and $P_2$ are associated with the compatible Poisson tensors
$\tilde P_1=e^{-t} P_1$ and $\tilde P_2=e^{-t}P_2$ on $\hat A$.
Since $\phi_0=\hat\d t$, compatibility between these Poisson
tensors implies that
$$
\brr{\tilde P_{1}^\sharp(\phi_0),\tilde P_2}_{\hat A}+\brr{\tilde
P_{2}^\sharp(\phi_0),\tilde P_1}_{\hat A}=0,
$$
or, using relation (\ref{relgauging1}),
$$
\brr{P_{1}^\sharp(\phi_0),P_2}^{\phi_0}+\brr{P_{2}^\sharp(\phi_0),P_1}^{\phi_0}=0.
$$
Now notice that
$$
\brr{P_{1}^\sharp(\phi_0),P_2}^{\phi_0}=\brr{P_{1}^\sharp(\phi_0),P_2}-i_{\phi_0}P_{1}^\sharp(\phi_0)\wedge
P_2=\brr{P_{1}^\sharp(\phi_0),P_2},
$$
so, compatibility between Jacobi bivectors also implies that
\begin{equation}
\brr{P_{1}^\sharp(\phi_0),P_2}+\brr{P_{2}^\sharp(\phi_0),P_1}=0.
\end{equation}

Now, let $(P^1_M=\rho^2 P_1,E^1_M=\rho (P_{1}^\sharp(\phi_0))$ and
$(P^2_M=\rho^2 P_2,E^2_M=\rho (P_{2}^\sharp(\phi_0))$ be the
Jacobi structures on $M$ induced by the triangular Jacobi
algebroids $(A,\phi_0,P_1)$ and $(A,\phi_0,P_2)$ (see
(\ref{jacobibivector:manifold:one}) and
(\ref{jacobibivector:manifold:two})).

Since $\rho$ is a Lie algebroid morphism, we have
\begin{align*}
    \brr{P^1_M,P^2_M}&=\brr{\rho^2P_1, \rho^2 P_2}=\rho^3\brr{P_1,P_2}=\rho^3(-P_{1}^\sharp(\phi_0)\wedge P_2-P_{2}^\sharp(\phi_0)\wedge
    P_1)\\
    &=-\rho(P_{1}^\sharp(\phi_0))\wedge \rho^2 P_{2}-\rho(P_{2}^\sharp(\phi_0))\wedge \rho^2
    P_1\\
    &=-E^1_M\wedge P^2_M-E_M^2\wedge P^1_M
\end{align*}
and
\begin{align*}
    \brr{E_M^1,P^2_M}+\brr{E_M^2,P^1_M}=\rho^2(\brr{P_{1}^\sharp(\phi_0),P_2}+\brr{P_{2}^\sharp(\phi_0),P_1})=0.
\end{align*}

So the given Jacobi structures on $M$ are compatible.
\end{proof}
%%%%%%%%%%%%%%%%%%%%%%%%%%%%%%%%%%% Jacobi-Nijenhuis algebroids %%%%%%%%%%%%%%%%%%%%%%%%%%%%%%%

\section{Jacobi-Nijenhuis algebroids}

We begin this section exposing  some well known results about
Nijenhuis operators and compatible Poisson structures on Lie
algebroids.

\subsection{Poisson-Nijenhuis Lie algebroids}
Let $(A,\br,\rho)$ be a Lie algebroid over a manifold $M$. Recall
that a \emph{Nijenhuis operator} is a bundle map $N:A\ra A$ (over
the identity) such that the induced map on the sections (denoted by
the same symbol $N$) has vanishing torsion:
\begin{equation}
\label{eq:Nijenhuis} T_N(X,Y):=[NX,NY]-N[X,Y]_N=0, \quad  X,Y\in
\X^1(A),
\end{equation}
where $\br_N$ is defined by
$$
[X,Y]_N:=[NX,Y]+[X,NY]-N[X,Y],\quad X,Y\in \X^1(A).
$$
Let us set $\rho_N:=\rho\circ N$. For a Nijenhuis operator $N$, one
easily checks  that the triple $A_N=(A,\br_N,\rho_N)$ is a new Lie
algebroid, and then $N:A_N\ra A$ is a Lie algebroid morphism.

Since $N$ is a Lie algebroid morphism, its transpose gives a chain
map of  com\-ple\-xes of dif\-fe\-ren\-tial forms
$N^\ast:(\Omega^k(A),\d_A)\ra(\Omega^k(A_N),\d_{A_N})$. Hence we
also have a map at the level of algebroid cohomology
$N^\ast:H^\bullet(A)\ra H^\bullet(A_N)$.

When the Lie algebroid $A$ is
 equipped with a Poisson structure $P$ and a
Nijenhuis operator $N$ which are \emph{compatible}, it is called a
\emph{Poisson-Nijenhuis Lie algebroid}.

The compatibility condition between $N$ and $P$ means that $NP$ is a
bivector field and
\[ \br_{N P}=\br_{P}^N,\]
where $\br_{NP}$ is the  bracket defined by the bivector field
$NP\in\X^2(A)$, and $\br_{P}^N$ is the bracket obtained from the
Lie bracket $\br_{N}$ by the Poisson bivector $P$.

As a consequence, $NP$ defines a new Poisson structure on $A$,
compatible with $P$:
\[ [P,NP]=[NP,NP]=0,\]
and one has a commutative diagram of morphisms of Lie algebroids:
\[
\xymatrix{
(A^\ast,[\cdot,\cdot]_{NP})\ar[rr]^{{N}^*}\ar[dd]_{P^\sharp}\ar[ddrr]^{NP^\sharp}
&&
(A^\ast,[\cdot,\cdot]_{P})\ar[dd]^{P^\sharp}\\
\\
(A,[\cdot,\cdot]_{N})\ar[rr]^{{N}}&&(A,[\cdot,\cdot]) }
\]

In fact, we have a whole hierarchy ${N^kP}$ ($k\in\N$) of pairwise
compatible Poisson structures on $A$.

\subsection{Jacobi-Nijenhuis algebroids}

Let $(A,{\phi_0})$ be a Jacobi algebroid and $N$ a Nijenhuis
operator on $A$. The definition of the Lie algebroid structure on
$\hat A=A\times \R$ given by (\ref{anchorhatA}) and
(\ref{brackethatA}) allows us to say that $N$ is also a Nijenhuis
operator on $\hat A$. So we have an additional Lie algebroid
structure on $\hat A$, $\hat A_N$.

\begin{prop}
The 1-form $\phi_1=N^\ast {\phi_0}$ is a 1-cocycle of $A_N$. The Lie
algebroid structure $\hat A_N$ coincides with the Lie algebroid
structure  on $\hat A$ induced from $A_N$  by $\phi_1$.
\end{prop}

\begin{proof}
First notice that, since $N: A_N\ra A$ is a Lie algebroid morphism,
 ${\ds \d_{N} \phi_1=\d_{N} N^\ast {\phi_0}=N^\ast (\d
{\phi_0})=0}$, and then  $\phi_1$ is a 1-cocycle of $A_N$. Besides,
for $X, Y\in\X^1(A)$, we have $NX,\, NY\in\X^1(A)$,
\begin{align*}
    \brr{X,Y}_{\hat A_ N}&=\brr{NX,Y}_{\hat A}+\brr{X,NY}_{\hat A}-N\brr{X,Y}_{\hat A}\\
    &= \brr{NX,Y}+\brr{X,NY}-N\brr{X,Y}=\brr{X,Y}_N
\end{align*}
and
\begin{align*}
    \hat \rho_N(X)&=\hat\rho\circ N(X)=\rho(NX)+\lan {\phi_0},NX\ran \frac{\partial}{\partial t}\\
                  &= \rho\circ N(X)+\lan N^\ast{\phi_0}, X\ran \frac{\partial}{\partial t}\\
                  &= \widehat{\rho_N}(X).
\end{align*}
Since $\X(\hat A)$, as  $C^\infty(M\times \R)$-module, is generated
by $\X(A)$, we conclude that $\hat A_N$ and the Lie algebroid
structure on $\hat A$ induced   from $A_N$ by $\phi_1$ are the same.
\end{proof}

In fact we have a whole sequence of Lie algebroid structures on
$\hat A$ given by $N^k$ or, equivalently, by the 1-cocycle of
$A_{N^k}$, $\phi_k=N^{\ast\,k}{\phi_0}$:
\begin{equation}\label{hierachyLiestructures}
 \hat A_{N^k}=(\hat A,\, \br_{N^k}, \hat\rho_{N^k}=\hat \rho\circ
N^k),\quad k\in\N.
\end{equation}

\

Now suppose $P\in\X^2(A)$ is a Jacobi bivector, i.e., a bivector
field such that $ \brr{P,P}^{\phi_0}=0. $ If $NP$ is a bivector on
$A$, we can consider the bracket on $A^\ast$ obtained from
$(A,{\phi_0})$ by  $NP$:

\begin{equation}\label{bracket1}
\brr{\alpha,\beta}_{NP}=\lie^{{\phi_0}}_{NP^\sharp{\al}}
\beta-\lie^{{\phi_0}}_{NP^\sharp{\beta}}
    \al-\d^{{\phi_0}}NP(\al,\beta), \,\,\,\, \al,\beta \in \X^1(
    A^\ast).
\end{equation}

On the other hand, we can also consider the bracket on $A^\ast$
obtained  from $(A_N,\phi_1=N^\ast{\phi_0})$ by the Jacobi bivector
$P$:

\begin{equation}\label{bracket2}
    \brr{\al,\beta}_P^N=\lie^{N,\,\phi_1}_{P^\sharp{\al}} \beta-\lie^{N,\,\phi_1}_{P^\sharp{\beta}}
    \al-\d_N^{\phi_1}P(\al,\beta), \,\,\,\, \al,\beta \in \X^1(
    A^\ast),
\end{equation}
where $\lie^{N,\,\phi_1}$ is the $\phi_1$- Lie derivative on $A_N$.

\begin{defi}\label{compatible:PandN}
The Jacobi bivector $P$ and the Nijenhuis operator $N$ are
\textbf{compatible} if the  following two conditions are satisfied:
\begin{enumerate}
  \item ${\ds NP=PN^\ast}$;
  \item the brackets $\br{}_{NP}$ and $\br{}_P^N$, given by (\ref{bracket1}) and (\ref{bracket2}), coincide.
\end{enumerate}
In this case, the Jacobi algebroid $(A,{\phi_0})$ is said to be a
\textbf{Jacobi-Nijenhuis algebroid} and is denoted by
$(A,{\phi_0},P,N)$.
\end{defi}

\

The compatibility between $N$ and $P$ can be expressed by the
vanishing of a suitable {\em concomitant}.

On a Jacobi algebroid $(A,{\phi_0})$ consider  a Nijenhuis operator
$N$ and a Jacobi bivector $P$ such that $NP$ is a bivector.
Following \cite{KosmannMagri}, we define the \emph{concomitant of
$P$ and $N$} as
\begin{align}\label{concomitant}
    C(P,N)(\al,\beta)=\brr{\alpha,\beta}_{NP}-\brr{\alpha,\beta}_P^N, \quad
    \al,\beta\in\Omega^1(A),
\end{align}
where $\br_{NP}$ and $\br_P^N$ are the brackets on $A^\ast$ given by
(\ref{bracket1}) and (\ref{bracket2}), respectively. We immediately
see that Condition (2) on Definition \ref{compatible:PandN} is
equivalent to $C(P,N)=0$.

 A direct
computation gives the following equalities, with $\al, \beta \in
\Omega^1(A)$:
\begin{align*}
\brr{\alpha,\beta}_{NP}= e^t\brr{\alpha,\beta}_{N\tilde P}-\lan
\alpha, NP^\sharp({\phi_0})\ran \beta+\lan \beta,
NP^\sharp({\phi_0})\ran \alpha,
\end{align*}
and
\begin{align*}
    \brr{\al,\beta}_P^N= e^t\brr{\alpha,\beta}_{\tilde P}^N-\lan \alpha, P^\sharp(\phi_1)\ran \beta+\lan \beta,
P^\sharp(\phi_1)\ran \alpha,
\end{align*}
where $\br_{\tilde P}^N$ is the bracket on $\hat A^\ast$ obtained
from $\hat A_N$ by $\tilde P$ and $\br_{N \tilde{P}}$ is the bracket
on $\hat A^\ast$ obtained from $\hat A$ by the bivector $N \tilde
P$.

Recall that compatibility  between the Poisson bivector $\tilde P$
and the Nijenhuis operator $N$, on the Lie algebroid $\hat{A}$,
means that $N\tilde P$ is a bivector and $\hat{C}(\tilde P,N)=0$,
where $\hat{C}(\tilde P,N)$ is the concomitant of $\tilde{P}$ and
$N$. Observing that
\begin{align*}
    C(P,N)(\al,\beta)=e^t \hat{C}(\tilde P, N)(\al,\beta),\quad
    \al,\beta\in\Omega^1(A),
\end{align*}
and also that
\begin{align*}
    C(P,N)(\phi_0,\al)=e^t \hat{C}(\tilde P, N)(\hat{d}t,\al),\quad
    \al \in\Omega^1(A),
\end{align*}
we conclude
%, from the $C^\infty(M)$-linearity of $C$ and the
%$C^\infty(M \times \R)$-linearity of $\hat{C}$,
 that $N$ and $P$
are compatible (on $(A, \phi_0)$) {\em if and only if}
 $N$ and $\tilde{P}$ are compatible (on $\hat{A}$).

 The Poisson-compatibility between $N$ and $\tilde{P}$ implies
$\brr{N \tilde{P}, \tilde{P}}_{\hat{A}}=0$ (see
\cite{KosmannMagri}). From $\brr{NP,P}^{\phi_0}=e^{3t}\brr{N
\tilde{P}, \tilde{P}}_{\hat{A}}$, we conclude that the
Jacobi-compatibility between $N$ and $P$ implies
$\brr{NP,P}^{\phi_0}=0$.

\

\begin{prop}
On a Jacobi-Nijenhuis algebroid $(A,{\phi_0},P,N)$, we have a
hierarchy of compatible Jacobi bivectors.
\end{prop}

\begin{proof}
Consider the Poisson bivector $\tilde P=e^{-t}P$ on $\hat A$ . As
we have already seen, the Jacobi-com\-pa\-ti\-bi\-li\-ty  between
$N$ and $P$ is equivalent
 to Poisson-compatibility between $N$ and
$\tilde P$ and we have a hierarchy of compatible Poisson bivectors
$N^k\tilde P$, $k\in \N$, on $\hat A$. This hierarchy induces a
hierarchy of compatible Jacobi bivectors on $A$, $N^kP$:
\begin{equation}\label{compcondbivectores}
\brr{N^iP, N^jP}^{{\phi_0}}=0, \quad (i,j\in\N).
\end{equation}
\end{proof}

\begin{corol}\label{Prop:hier:jacobi:M} The  Jacobi-Nijenhuis algebroid $(A,{\phi_0}, P, N)$
defines a
 hierarchy of compatible Jacobi structures
on $M$.
\end{corol}

\begin{proof}
This is an immediate consequence of the above proposition and
theorem \ref{compatible:jacobi:bivectors}.
\end{proof}

Also the compatibility conditions define a sequence of Lie algebroid
structures on $A^\ast$.

\begin{theorem}
Let $(A,{\phi_0},P,N)$ be a Jacobi-Nijenhuis algebroid. Then
$A^\ast$ has a hierarchy of Jacobi algebroid structures $(A^\ast,
X_k)$, such that $((A,\phi_i), (A^\ast,X_k))$,
$\phi_i=N^{\ast\,i}\phi_0$ and $X_k=N^kX_0$, $i\leq k$, $k\in\N$,
are triangular Jacobi bialgebroids.
\end{theorem}

\begin{proof}
Last proposition guarantees that $N^iP$, $i\in\N$, is a hierarchy of
compatible Jacobi bivectors.

Each one of the Poisson bivectors $N^k\tilde
P=\widetilde{N^kP}=e^{-t}N^kP$ defines a Lie algebroid structure on
$\hat A^\ast$,
\begin{equation*}
  \hat A^\ast_{N^k}=(\hat A^\ast, \br_{N^k\tilde P}, \hat\rho_{k\,\ast}=\hat\rho\circ
  N^k\tilde
  P^\sharp),
\end{equation*}
and a Lie algebroid structure on $A^\ast$,
\begin{equation*}
  A^\ast_{N^kP}=( A^\ast, \br_{N^kP}, \rho_{k\,\ast}=\rho\circ N^k
  P^\sharp),
\end{equation*}
where
\begin{equation*}
    \brr{\alpha,\beta}_{N^kP}=e^{-t}\brr{e^t\alpha, e^t\beta}_{N^k\tilde
    P}, \quad \alpha,\beta\in \X^1(A^\ast).
\end{equation*}
Each Lie algebroid structure $A^\ast_{N^kP}$ coincides with the Lie
algebroid structure obtained from the Jacobi algebroid
$(A,\phi_{k-i})$ by the Jacobi bivector $N^iP$, $i=1,\ldots,k$. So,
the pairs $((A,\phi_{k-i}), (A^\ast,X_k))$, $i=1,\ldots,k$ are
triangular Jacobi bialgebroids.
\end{proof}

As in the Poisson case, $N^\ast$ is a Nijenhuis operator of $A^\ast$
and we have a
 commutative relation between duality by $P$ and deformation along $N^\ast$.

\begin{prop} Let $(A, {\phi_0}, P,N)$ be a Jacobi-Nijenhuis algebroid and
consider the Lie algebroid structure on $A^\ast$ given by
(\ref{bracketdualA}) and (\ref{anchordualA}).  The operator $N^\ast$
is a Nijenhuis operator on $A^\ast$.
\end{prop}

\begin{proof}
Since relation (\ref{relgauging3}) holds and $N^\ast$ is a Nijenhuis
operator on $\hat A^\ast$,  we have
\begin{align*}
    T^{A^\ast}_{N^\ast}(\alpha,&\beta)
    =[N^\ast\alpha,N^\ast\beta]_{P}-N^\ast([N^\ast\alpha,\beta]_{P}+[\alpha,N^\ast\beta]_{P}-N^\ast[\alpha,\beta]_{P})\\
    &=e^{-t}\left([N^\ast\hat\alpha,N^\ast\hat\beta]_{\tilde P}-N^\ast([N^\ast\hat \alpha,\hat \beta]_{\tilde P}
    +[\hat \alpha,N^\ast\hat \beta]_{\tilde P}
    -N^\ast[\hat \alpha,\hat \beta]_{\tilde P})\right)\\
    &=e^{-t}T^{\hat A^\ast}_{N^\ast}(\hat \alpha,\hat \beta)=0,
\end{align*}
where $\hat\alpha=e^t\al$,  $\hat\beta=e^t\beta$ and
$\al,\beta\in\Omega^1(A)$.
\end{proof}

This way $A^\ast$ can be deformed by $N^\ast$ into $A^\ast_{N^\ast}$
and one can easily check that this is exactly the Lie algebroid
$A^\ast_{NP}$.

\begin{prop}
The Lie algebroid $A^\ast_{N^kP}$ coincides  with the Lie algebroid
  $A^\ast_{N^{\ast\, k}}$, obtained from $A^\ast$ by deformation along $N^{\ast\, k}$.
\end{prop}

\

We finish this section, showing that the definition of strong (or
strict)Jacobi-Nijenhuis structure defined for Jacobi manifolds in
\cite{IglMarr2} and \cite{Cos} can be recovered in this framework.

%%%%%%%%%%%%%%%%%%%%%%%%%%%%%%%% STRONG JACOBI-NIJENHUIS MANIFOLD %%%%%%%%%%%%%%%%%%%%%%%%%%%%%%%%%

\begin{ex}Consider a Jacobi manifold $(M, (\Lambda,E))$ and the Lie
algebroid  $A=TM\times \R$  defined in section 2.2. A \emph{strong
(or strict) Jacobi-Nijenhuis structure} on $M$ is given by a
Nijenhuis operator on $A$, $\mathcal{N}$, compatible with
$(\Lambda,E)$ in the following sense:

(i) ${\ds
\mathcal{N}\circ(\Lambda,E)^\sharp=(\Lambda,E)^\sharp\circ\mathcal{N}^\ast}$.
This condition defines a new skew-symmetric bivector $\Lambda_1$
and a vector field $E_1$ such that ${\ds (\Lambda_1,
E_1)^\sharp=\mathcal{N}\circ(\Lambda,E)^\sharp}$.

(ii) The concomitant of $(\Lambda,E)$ and $\mathcal{N}$, ${\ds
\mathcal{C}((\Lambda,E),\mathcal{N})}$, identically vanishes.

The concomitant ${\ds \mathcal{C}((\Lambda,E),\mathcal{N})}$ is
given in \cite{IglMarr2, Cos} by
\begin{align*}
\mathcal{C}((\Lambda,E),\mathcal{N})((\al,f),(\beta,g))&=\brr{(\al,f),(\beta,g)}_{(\Lambda_1,E_1)}-\brr{\mathcal{N}^\ast(\al,f),(\beta,g)}_{(\Lambda,E)}\\
&-\brr{(\al,f),\mathcal{N}^\ast(\beta,g)}_{(\Lambda,E)}+\mathcal{N^\ast}\brr{(\al,f),(\beta,g)}_{(\Lambda,E)},
\end{align*}
for $(\al,f),(\beta,g)\in\Omega^1(M)\oplus C^\infty(M)$, where the
brackets $\br_{(\Lambda,E)}$ and $\br_{(\Lambda_1,E_1)}$ are
defined in (\ref{bracket:TdualMXR}).

The  concomitant can be rewritten as
\begin{align*}
\mathcal{C}((\Lambda,E),\mathcal{N})((\al,f),(\beta,g))&=\brr{(\al,f),(\beta,g)}^\mathcal{N}_{(\Lambda,E)}
-\brr{(\al,f),(\beta,g)}_{(\Lambda_1,E_1)}
\end{align*}
and we obtain the symmetric of (\ref{concomitant}).

We conclude  that a strong Jacobi-Nijenhuis structure is a pair of
compatible Nijenhuis and Jacobi structures in the sense of
 definition
\ref{compatible:PandN}.
\end{ex}
%%%%%%%%%%%%%%%%%%%%%%%%%%%%%%%%%%%%%%%%%%%%%%%%%%%%%%%%%%%%%%%%%%%%%%%%%%%

%
\section{Modular classes of Jacobi algebroids}

\subsection{Modular class of a Lie algebroid}
Let $(A, \br, \rho)$ be a Lie algebroid over the manifold $M$. For
simplicity we will assume that both $M$ and $A$ are orientable, so
that there exist non-vanishing sections
$\eta\in\X^{\mathrm{top}}(A)$ and $\mu\in\Omega^{\mathrm{top}}(M)$.

The \emph{modular form} of the Lie algebroid $A$ with respect to
$\eta\otimes\mu$ (see \cite{EvLuWe}) is the 1-form
$\xi_A^{\eta\otimes\mu}\in\Omega^1(A)$, defined by
\begin{equation}\label{defi:modular:form}
    \langle \xi_A^{\eta\otimes\mu},X\rangle \eta\otimes \mu=\lie_X\eta\otimes
    \mu+\eta \otimes\lie_{\rho(X)}\mu, \quad X\in\mathfrak{X}^1(A).
\end{equation}
This is a 1-cocycle of the Lie algebroid cohomology of $A$. If one
makes a different choice of sections $\eta'$ and $\mu'$, then
$\eta'\otimes\mu'=f\eta\otimes \mu$, for some non-vanishing smooth
function $f\in C^\infty(M)$. One checks easily that the modular form
$\xi_A^{\eta'\otimes\mu'}$ associated with this new choice is given
by:
\begin{equation}
\label{eq:change:section}
\xi_A^{\eta\otimes\mu}=\xi_A^{\eta'\otimes\mu'}-\d\ln|f|,
\end{equation}
so that the cohomology class $[\xi_A^{\eta\otimes\mu}]\in H^1(A)$ is
independent of the choice of $\eta$ and $\mu$. This cohomology class
is called the \emph{modular class of the Lie algebroid}  $A$ and we
will denoted it by $\m A:=[\xi_A^{\eta\otimes\mu}]$.

\subsection{Modular classes of a Jacobi algebroid}

Let $(A, {\phi_0})$ be a Jacobi algebroid of rank $n$. The
Schouten-Jacobi bracket $\br^{\phi_0}$, given by
(\ref{Jacobibracket}), allows us to define a representation of $A$
on $Q_A=\X^{n}(A)\otimes \Omega^{\mathrm{top}}(M)$.

\begin{prop}
Let $(A, {\phi_0})$  be a Jacobi algebroid. The bilinear map
$D^{\phi_0}: \X^1(A)\otimes Q_A \rightarrow Q_A$ defined by
\begin{equation}\label{phirepresentation}
 D^{\phi_0}_X\,(\eta\otimes\mu)=\brr{X,\eta}^{\phi_0}\otimes\mu + \eta\otimes
\lie_{\rho(X)}\mu,
\end{equation}
is a representation of the Lie algebroid $A$ on $Q_A$.
\end{prop}

\begin{proof}
By definition of $\br^{\phi_0}$, we have
$$
D^{\phi_0}_X\,(\eta\otimes\mu)=\left(\brr{X,\eta}-(n-1)\lan{\phi_0},X\ran\eta\right)\otimes\mu
+ \eta\otimes \lie_{\rho(X)}\mu,
$$
so $D^{\phi_0}=D-(n-1){\phi_0}$, where $D$ is the representation of
$A$ on $Q_A$ considered in \cite{EvLuWe} to define the modular class
of the Lie algebroid $A$.

Obviously, for $f\in C^\infty(M)$, $X,Y\in\X^1(A)$ and $s\in
\Gamma(Q_A)$, $D^{\phi_0}$ satisfies
\begin{equation*}
    D^{\phi_0}_{fX}s=fD^{\phi_0}_X s,
\end{equation*}
and
\begin{equation*}
    D^{\phi_0}_X(fs)=fD^{\phi_0}_Xs +
    (\rho(X)f)s.
\end{equation*}
Moreover, since $D$ is a representation and ${\phi_0}$ a 1-cocycle
of $A$,
\begin{align*}
    D^{\phi_0}_X(D^{\phi_0}_Y s)&-D^{\phi_0}_Y(D^{\phi_0}_X s)=\\
    &= D^{\phi_0}_X(D_Y s-(n-1)\lan {\phi_0}, Y\ran s) - D^{\phi_0}_Y(D_X s-(n-1)\lan {\phi_0}, X\ran
    s)\\
    &= (D_XD_Y-D_YD_X)s -(n-1)(\rho(X){\phi_0}(Y)-\rho(Y){\phi_0}(X))s\\
    &= D_{\brr{X,Y}}s-(n-1){\phi_0}(\brr{X,Y})s=D^{\phi_0}_{\brr{X,Y}}s.
\end{align*}
We conclude that $D^{\phi_0}$ is a representation of $A$ on $Q_A$.
\end{proof}

\begin{defi}
The \textbf{modular form of the Jacobi algebroid} $(A, {\phi_0})$
with respect to $\eta\otimes \mu$ is the $A$-form
$\xi_A^{{\phi_0},\,\eta\otimes \mu}$ defined by
$$
\xi_A^{{\phi_0},\,\eta\otimes \mu}=\xi_A^{\eta\otimes
\mu}-(n-1){\phi_0}.
$$
\end{defi}

Again, the cohomology class of a modular form is independent of the
section of $Q_A$ chosen.

\begin{defi}
 The
\textbf{modular class of the Jacobi algebroid} $(A,{\phi_0})$  is
the cohomology class of a modular form. It  will be denoted  by
$\m^{\phi_0} A=\brr{\xi_A^{{\phi_0},\,\eta\otimes\mu}}$.
\end{defi}

Obviously $\m^{\phi_0} A=\m A$ if and only if ${\phi_0}$ is exact.

\subsection{Relation between the modular classes of $A$ and $\hat
A$ and of $A^\ast$ and $\hat A^\ast$.}\label{sectionrelations}

Let $(A,{\phi_0})$ be a Jacobi algebroid of rank $n$. In this
section we compute modular forms of $\hat A$ and $\hat A^\ast$ (in
the triangular case) and we establish relations between them and the
modular forms of  $A$ and $A^\ast$.
 Let $\eta\in\X^{n}(A)$
and $\mu\in\Omega^{\mathrm{top}}(M)$, then $\eta$ is also a
$n$-section of $\hat A$  and $\tilde \mu=\mu\wedge d t$ is a
volume form of $M\times \R$.

The Lie bracket on $\hat A$ coincides with the Lie bracket on $A$
for time-independent multivectors, so
\begin{align*}
\xi_{\hat A}^{\eta\otimes\tilde\mu}(X)\eta\otimes\tilde\mu&=
\brr{X,\eta}_{\hat A}\otimes \tilde \mu+  \eta \otimes \lie_{\hat
\rho(X)}\tilde \mu\\
&=\brr{X,\eta}\otimes \tilde \mu+  \eta \otimes \lie_{\rho(X)+\lan
{\phi_0}, X\ran\frac{\p}{\p t}}\tilde\mu, \quad X\in\X^1(A).
\end{align*}
Since
\begin{align*}
\lie_{\rho(X)+\lan {\phi_0}, X\ran\frac{\p}{\p t}}(\mu\wedge d
t)=\lie_{\rho(X)}\mu\wedge d t,
\end{align*}
we have
\begin{align}
\xi_{\hat
A}^{\eta\otimes\tilde\mu}(X)\eta\otimes\tilde\mu&=\brr{X,\eta}\otimes
\tilde\mu+ \eta\otimes \lie_{\rho(X)}\mu\wedge
d t\nonumber \\
& = \xi_A^{\eta\otimes\mu}(X) \eta\otimes
\tilde\mu.\label{relation:modular:hatA}
\end{align}

Now consider the section of $\wedge^n\hat A$, $\tilde
\eta=e^{-(n-1)t}\eta$.  Using relation (\ref{relgauging1}) and
definition of the Schouten-Jacobi bracket $\br^{\phi_0}$, we find
that the modular form of $\hat A$ with respect to $\tilde
\eta\otimes\tilde \mu$ is given by
\begin{align*}
\xi_{\hat
A}^{\tilde\eta\otimes\tilde\mu}(X)\tilde\eta\otimes\tilde\mu&=
 \left(\xi_A^{\eta\otimes\mu}(X)-(n-1)\lan {\phi_0}, X
\ran\right)\tilde \eta\otimes
\tilde\mu\\&=\xi_A^{{\phi_0},\,\eta\otimes\mu}(X)\tilde \eta\otimes
\tilde\mu.
\end{align*}

\begin{prop}Let $(A, {\phi_0})$ be a Jacobi algebroid, then
$$
\brr{\xi_{\hat A}}=\brr{\xi_A^{\phi_0}}=\brr{\xi_A}.
$$
\end{prop}

It is clear that  the cohomology considered in the previous
proposition is the $\hat A$-cohomology. In $\hat A$ the 1-form
${\phi_0}$ is exact, ${\phi_0}=\hat \d t$ and, generally, this is
not the case in $A$.

Now suppose we also have a Jacobi bivector $P$ on $(A,{\phi_0})$. We
saw that it induces a Poisson structure on $\hat A$, a Lie algebroid
structure on $\hat A^\ast$ and another one on $A^\ast$. Consider
$\nu\in\X^{n}(A^\ast)$ a top-section on $A^\ast$ and $\mu$ a volume
form on $M$.

\begin{prop}
The modular form of the Lie algebroid $A^\ast$ with respect to
$\nu\otimes \mu$ is given by
\begin{align}\label{relacao:modular:dualAum}
\xi_{A^\ast}^{\nu\otimes\mu}(\alpha)=e^t \xi_{\hat
A^\ast}^{\hat\nu\otimes\tilde \mu}(\alpha),\quad
\alpha\in\X^1(A^\ast),
\end{align}
with $\hat\nu=e^{nt}\nu$ and $\tilde\mu=\mu\wedge d t$.

The modular form of the Jacobi algebroid $(A^\ast, X_0)$, where
$X_0=-P^\sharp({\phi_0})$, with respect to $\nu\otimes \mu$ is given
by
\begin{equation}\label{relacao:modular:dualA}
\xi_{A^\ast}^{X_0,\,\nu\otimes\mu}= e^t\xi_{\hat
A^\ast}^{\nu\otimes\tilde\mu}+X_0.
\end{equation}
\end{prop}

\begin{proof}
By definition of modular form and relation (\ref{relgauging3}), for
$\alpha\in\Omega^1(A)$, we have
\begin{align*}
   e^t\xi_{\hat A^\ast}^{\hat \nu\otimes \tilde\mu}( \alpha)\hat\nu\otimes \tilde\mu
    &=\xi_{\hat A^\ast}^{\hat \nu\otimes \tilde\mu}(\hat \alpha)\hat\nu\otimes \tilde\mu \\&=
                   \brr{e^t\alpha,e^{nt}\nu}_{\tilde P}\otimes
                   \tilde\mu + \hat\nu \otimes \lie_{\hat
                   \rho_\ast(e^t\alpha)}(\mu\wedge d t)\\
    &= e^{nt}{\brr{\alpha,\nu}_{P}}\otimes \tilde\mu +
    \hat\nu \otimes \lie_{\hat\rho(P^\sharp(\alpha))}(\mu\wedge d t)\\
    &=e^{nt}{\brr{\alpha,\nu}_{P}}\otimes \tilde\mu +
    \hat\nu \otimes \lie_{\rho(P^\sharp(\alpha))}(\mu)\wedge d t\\
    &= \xi_{A^\ast}^{\nu\otimes\mu}(\alpha) (\hat
    \nu\otimes\tilde\mu).
\end{align*}
So
$$
\xi_{A^\ast}^{\nu\otimes\mu}(\alpha)= e^t\xi_{\hat A^\ast}^{ \hat
\nu\otimes \tilde\mu}(\alpha), \quad \alpha\in\Omega^1(A).
$$

\

Since $\nu$ is a $n$-form of $A$, we have $\al\wedge
i_{X_0}\nu=\lan\al,X_0\ran\nu$, $\al\in\Omega^1(A)$, and using
relation (\ref{bracketchapeuX0})  we obtain
\begin{align*}
\brr{\alpha,\nu}_{\tilde P}&=e^{-t}\left(\brr{\al,\nu}_P-n\lan \al,X_0\ran \nu+\al\wedge i_{X_0}\nu\right) \\
&=e^{-t}\left(\brr{\al,\nu}_P-(n-1)\lan \al,X_0\ran \nu\right)\\
&=e^{-t}\brr{\alpha,\nu}_{P}^{X_0}.
\end{align*}

Also we have  \begin{align}
\lie_{\hat\rho_\ast(\alpha)}\tilde\mu&=\lie_{e^{-t}\hat\rho
(P^\sharp\al)}\tilde \mu
=e^{-t}\lie_{\hat\rho(P^\sharp\al)}\tilde \mu + \lan \hat\d e^{-t}, \hat\rho(P^\sharp\al)\ran \tilde \mu\nonumber\\
&= e^{-t} \left(\lie_{\rho(P^\sharp \al)+\lan
\phi_0,P^\sharp\al\ran\frac{\p}{\p t}}\tilde\mu- \lan
\phi_0,P^\sharp\al\ran\tilde\mu\right)\nonumber\\
&=e^{-t} \left(\lie_{\rho(P^\sharp \al)}\tilde \mu+\lie_{\lan
\phi_0,P^\sharp\al\ran \frac{\p}{\p t}}\tilde\mu- \lan
\phi_0,P^\sharp\al\ran\tilde\mu\right)\nonumber\\
&= e^{-t}(\lie_{\rho(P^\sharp \alpha)}\mu\wedge d t-\lan
\alpha,X_0\ran)\tilde\mu\nonumber\\
&=e^{-t}\left(\lie_{\rho(P^\sharp(\alpha))}\mu\wedge d t -\lan
\al, X_0\ran \tilde\mu\right).\label{liederivative:anchored}
\end{align}
These relations imply that
\begin{align*}
\xi_{\hat A^\ast}^{\nu\otimes\tilde\mu}(\al)\nu\otimes\tilde\mu&
=\brr{\alpha,\nu}_{\tilde P}\otimes
                   \tilde\mu + \nu \otimes \lie_{\hat
                   \rho_\ast(\alpha)}(\mu\wedge d t)\\
&=e^{-t}\brr{\alpha,\nu}_{P}^{X_0}\otimes\tilde\mu+e^{-t}\nu\otimes\left(\lie_{\rho(P^\sharp(\alpha))}\mu\wedge
d
t -\lan \al, X_0\ran \tilde\mu\right)\\
&= e^{-t}\left( \xi_{A^\ast}^{X_0,\,\nu\otimes\mu}(\alpha)- \lan
\al, X_0\ran \right)\nu\otimes\tilde\mu
\end{align*}
and  relation (\ref{relacao:modular:dualA}) follows.
\end{proof}

%%%%%%%%%%%%%%%%%%%%%%% Relation with the modular vector field defined by Marrero %%%%%%%%%%%%%%%%%

\subsection{Relation with the modular vector field of a triangular Jacobi
bialgebroid}

 The definition of modular class of a
triangular Jacobi bialgebroid was given in \cite{IglLopMarPad}. In
this section we will present this definition using the approach we
have chosen, relating it with the modular field of the triangular
bialgebroid associated with the Jacobi bialgebroid.

Let $(\hat A, \tilde P )$ be the triangular Lie bialgebroid
associated with the triangular Jacobi algebroid $(A,{\phi_0}, P)$ of
rank  $n$  and
 $\nu$  a section of $\wedge^{n}A^\ast$. The modular  field of the
triangular Lie bialgebroid $(\hat A,\tilde P)$  with respect to
$\hat \nu = e^{nt}\nu$ (see \cite{Kosmann2}) is the section $\hat
X^{\hat\nu}$ of $\hat A$ given by
\begin{align*}
\hat X^{\hat\nu}(\alpha)\hat\nu&=-\alpha\wedge \hat\d\, i_{\tilde
P}\hat\nu=-\alpha\wedge \hat\d i_{e^{-t}P}(e^{nt}\nu)
\\
&= -\alpha\wedge \hat\d (e^{(n-1)t}i_P\nu) \\
&= -e^{(n-1)t}\alpha\wedge ((n-1) {\phi_0}\wedge i_P\nu + \d
i_P\nu), \quad \al\in \Omega^1(A).
\end{align*}

Comparing with the definition of $\mathcal{M}_{(A,{\phi_0},P)}^\nu$,
the modular vector field of the triangular Jacobi bialgebroid
$(A,{\phi_0},P)$ given in \cite{IglLopMarPad}, we notice that
\begin{equation}
\hat X^{\hat\nu}=e^{-t}\mathcal{M}_{(A,{\phi_0},P)}^\nu.
\end{equation}

Since $\br_{\tilde P}$ is generated by
$\partial_{\tilde{P}}=\hat{\d}i_{\tilde P}-i_{\tilde P}\hat\d$, we
have
\begin{equation}\label{modular:generator}
\hat X^{\hat\nu}(\al)\hat\nu=\brr{\alpha,\hat\nu}_{\tilde
P}+e^{-t}(i_P\d\alpha)\hat\nu.
\end{equation}
Moreover (see (\ref{liederivative:anchored})),
\begin{align*}
\lie_{\hat\rho(\tilde P)^\sharp\al}\tilde\mu&=
e^{-t}(\mathrm{div}_\mu\,\rho (P^\sharp\alpha)-\lan
\alpha,X_0\ran)\tilde\mu,
\end{align*}
where $\tilde \mu=\mu\wedge dt$, $\mu\in\Omega^{\mathrm{top}}(M)$.
Using the definition of modular form of a Lie algebroid
(\ref{defi:modular:form}):
$$
\xi_{\hat
A^\ast}^{\hat\nu\otimes\tilde\mu}(\alpha)\hat\nu\otimes\tilde\mu=\brr{\alpha,\hat
\nu}_{\tilde P}\otimes \tilde\mu + \hat\nu\otimes
\lie_{\hat\rho(\tilde P^\sharp(\al))}\tilde\mu
$$
and relation
 (\ref{relacao:modular:dualAum}),
we obtain
\begin{align}\label{relation:modular:com:Marrero}
    \xi_{A^\ast}^{\nu\otimes\mu}(\alpha)=\mathcal{M}_{(A,{\phi_0},P)}^\nu(\alpha)-i_P\d\al-\lan\alpha,X_0\ran+\mathrm{div}_\mu(\rho(P^\sharp(\al))).
\end{align}

On the other hand, notice that relation (\ref{modular:generator})
implies $\lan\hat \d f,\hat X^{\hat\nu}\ran\hat\nu=\brr{\hat\d
f,\hat\nu}_{\tilde P}$, $f\in C^\infty(M\times \R)$, so
\begin{align}\label{relation:modtrianbialgebroid:modLie}
    \hat\rho(\xi_{\hat A^\ast}^{\hat\nu\otimes\tilde\mu})=\hat\rho
    (\hat X^{\hat\nu})+ X^{T(M\times\R)},
\end{align}
where $X^{T(M\times\R)}$ is the modular vector field of the Poisson
manifold $M\times\R$ (endowed with the Poisson bivector induced from
the triangular Lie bialgebroid $(\hat A, \tilde P)$).

Since the 1-form ${\phi_0}$  is closed, we have
$$
\xi_{A^\ast}^{\nu\otimes\mu}(\phi_0)=\mathcal{M}_{(A,{\phi_0},P)}^\nu(\phi_0)-\mathrm{div}_\mu(\rho(X_0)),
$$
so
\begin{align*}
\hat\rho(\xi_{\hat
A^\ast}^{\hat\nu\otimes\tilde\mu})&=e^{-t}\hat\rho(\xi_{
A^\ast}^{\nu\otimes\mu})=e^{-t}\left(\rho(\xi_{
A^\ast}^{\nu\otimes\mu})+\lan \phi_0,
\xi_{A^\ast}^{\nu\otimes\mu}\ran \frac{\p}{\p t}\right)\\
&=e^{-t}\left(\rho(\xi_{ A^\ast}^{\nu\otimes\mu})
+\left(\mathcal{M}_{(A,{\phi_0},P)}^\nu(\phi_0)-\mathrm{div}_\mu
\rho(X_0)\right )\frac{\p}{\p t}\right).
\end{align*}
On another hand, ${\ds \hat\rho(\hat X^{\hat
\nu})=e^{-t}\left(\rho(\mathcal{M}_{(A,{\phi_0},P)}^\nu)+\mathcal{M}_{(A,{\phi_0},P)}^\nu(\phi_0)\frac{\p}{\p
t}\right)}$ and equation (\ref{relation:modtrianbialgebroid:modLie})
can be rewritten as
\begin{align}\label{relation:modtrianbialgebroid:modLie:dois}
   \rho(\xi_{A^\ast}^{\nu\otimes\mu})=\rho (\mathcal{M}^\nu_{(A,{\phi_0},P)})+
    e^t X^{T(M\times\R)}+ \mathrm{div}_{\mu} \rho(X_0)\frac{\p}{\p
    t}.
\end{align}

\

Let $(P_M,E_M)$ be the Jacobi structure on $M$ induced by the
triangular Jacobi algebroid $(A,\phi_0, P)$, i.e.,
$$
P_M(df,dg)=P(\d f,\d g),\quad E_M=\rho\circ P^\sharp({\phi_0}).
$$
The modular field of the Jacobi manifold $(M,P_M,E_M)$,
$V^{(P_M,E_M)}$, was introduced in \cite{Vaisman} and is defined as
$$
V^{(P_M,E_M)}=e^t X^{T(M\times\R)}.
$$
So,  equation (\ref{relation:modtrianbialgebroid:modLie:dois}) is
equivalent to
\begin{align*}%\label{relation:modtrianbialgebroid:modLie:tres}
    \rho(\xi_{A^\ast}^{\nu\otimes\mu})=\rho (\mathcal{M}_{(A,{\phi_0},P)})+
    V^{(P_M, E_M)}+ \mathrm{div}_{\mu} \rho(X_0)\frac{\p}{\p
    t}.
\end{align*}

%%%%%%%%%%%%%%%%%%%%%%%%%%%%%%%%%%% DUALITY OF MODULAR CLASSES %%%%%%%%%%%%%%%%%%%%%%%%%%%%%%%%%%

\subsection{Duality between modular classes of $A$ and  $A^\ast$}

Following the philosophy of this paper, we will find a relation
between the modular classes of the Jacobi algebroids $(A, {\phi_0})$
and $(A^\ast,X_0)$ using relations on the associated Lie
bialgebroid. So we begin by presenting some results about duality of
modular classes on Lie bialgebroids.

\begin{prop}
Let $(A,\br,\rho)$ be a Lie algebroid equipped with a Poisson
bivector $P$,  $(\br_P, \rho_\ast=\rho\circ P^\sharp)$ the Lie
algebroid structure induced by $P$ on $A^\ast$  and $\nu$ a
top-section on $A^\ast$. For all $\al\in\Omega^1(A)$, we have
\begin{align*}
\lie_{P^\sharp \al} \nu&=\brr{\al,\nu}_P + 2i_P (\d\alpha)\,\nu\\
&= -\brr{\al,\nu}_P-2\al\wedge \d i_P\nu.
\end{align*}
\end{prop}

\begin{proof}
Since $\nu$ is a top-section of $A^\ast$, using Cartan's formula, we
have
\begin{align}\label{relacao:um}
    \lie_{P^\sharp\alpha}\nu=\d i_{P^\sharp \al}\nu.
\end{align}
But $\al\wedge\nu=0$ and $
   i_P(\al\wedge\nu) = -i_{P^\sharp\al}\nu+\al\wedge i_P\nu
$, so $ i_{P^\sharp\al}\nu=\alpha\wedge i_P\nu. $ Substituting in
(\ref{relacao:um}) we have
\begin{align*}
    \lie_{P^\sharp\alpha}\nu=\d \alpha\wedge i_P\nu-\alpha\wedge
    \d i_P\nu.
\end{align*}

Again because $\nu$ is a top-section, we have that
$i_P(\d\alpha\wedge\nu)=0$, so $i_P(\d\al)\nu=\d\alpha\wedge i_P\nu$
and
\begin{equation}\label{relacao:dois}
\lie_{P^\sharp\alpha}\nu=i_P(\d\al)\nu-\al\wedge\d i_P\nu.
\end{equation}

 On the other hand, using the fact that $\partial_P=\brr{\d,i_P}$ is a
 generator of the Gerstenhaber algebra of $A^\ast$, we have
\begin{align*}
    \brr{\al,\nu}_P&= -i_P (\d\al)\,\nu-\al\wedge\d i_P\nu\\
    &=\lie_{P^\sharp\al}\nu-2 i_P(\d\alpha)\nu
\end{align*}
or, equivalently,  ${\ds
\brr{\al,\nu}_P=-\lie_{P^\sharp\al}\nu-2\alpha\wedge\d i_P\nu }$.
\end{proof}

\begin{prop}
Let $(A,A^\ast,P)$ be a triangular Lie bialgebroid. Then
\begin{equation}\label{relacao:tres}
    P^\sharp\xi_{A}^{\eta\otimes\mu}(\alpha)=-\xi_{A^\ast}^{\nu\otimes\mu}(\al)-2\lan\alpha\wedge\d
i_P\nu,\eta\ran,\quad \al\in\Omega^1(A),
\end{equation}
where $\mu$ is a volume form of $M$, $\eta\in\X^{\mathrm{top}}(A)$
and  $\nu\in\Omega^{\mathrm{top}}(A)$ such that $\lan
\nu,\eta\ran=1.$
\end{prop}

\begin{proof}
Since $\lan \nu,\eta\ran=1$, we have
$$
\lan \nu, \brr{X,\eta}\ran=-\lan\lie_X\nu, \eta\ran, \quad
X\in\mathfrak{X}^1(A),
$$
and
\begin{align*}
    \xi_A^{\eta\otimes\mu}(P^\sharp\al)\eta\otimes\mu&=\brr{P^\sharp\al,
    \eta}\otimes\mu + \eta\otimes \lie_{\rho(P^\sharp\al)}\mu\\
    &= \lan \nu, \brr{P^\sharp\al,\eta}\ran \eta\otimes\mu + \eta\otimes
    \lie_{\rho(P^\sharp\al)}\mu\\
    &=-\lan\lie_{P^\sharp\al}\nu, \eta\ran{\eta\otimes\mu}+ \eta\otimes
    \lie_{\rho(P^\sharp\al)}\mu\\
    &= \lan \brr{\al,\nu}_P+2\,\al\wedge \d i_P\nu, \eta\ran{\eta\otimes\mu} + \eta\otimes
    \lie_{\rho(P^\sharp\al)}\mu\\
    &= (\xi_{A^\ast}^{\nu\otimes\mu}(\al) +2\,\lan\al\wedge \d i_P\nu,
    \eta\ran)\eta\otimes\mu.
\end{align*}
So,  ${\ds
P^\sharp(\xi_A^{\eta\otimes\mu})(\al)=-\xi_{A^\ast}^{\nu\otimes\eta}(\al)-2\lan\alpha\wedge\d
i_P\nu,\eta\ran}$.
\end{proof}

 \

Now let $(A,{\phi_0})$ be a Jacobi algebroid of rank $n$ and $P$ a
Jacobi bivector on $A$. The pair $(\hat A,\tilde P)$ is a triangular
Lie bialgebroid  and we can use the previous proposition to relate
the modular classes of $\hat A$ and $\hat A^\ast$.

Consider $\eta\in\X^n(A)$ and $\nu\in\Omega^n(A)$ such that
$\lan\nu,\eta\ran=1$, then we have
\begin{align*}
    \tilde P^\sharp(\xi_{\hat A}^{\eta\otimes\tilde\mu})(\al)= -\xi_{\hat
    A^\ast}^{\nu\otimes\tilde\mu}(\al)-2\lan\al \wedge\hat\d i_{\tilde P}\nu,\eta\ran.
\end{align*}
Relations (\ref{relation:modular:hatA}) and
(\ref{relacao:modular:dualA}) imply that
\begin{align*}
    \tilde P^\sharp(\xi_{A}^{\eta\otimes\mu})(\al)= -e^{-t}(\xi_{
    A^\ast}^{X_0,\,\nu\otimes\mu}(\al)-\lan\al,X_0\ran) - 2\lan\al \wedge\hat\d i_{\tilde P}\nu,\eta\ran,
\end{align*}
and, since $\alpha\wedge\hat\d i_{\tilde P}\nu=\alpha\wedge \hat\d
(e^{-t}i_P\nu)=e^{-t} \left(\alpha\wedge\d i_P\nu -
P({\phi_0},\alpha)\nu\right)$, we have
\begin{align*}
 P^\sharp(\xi_{A}^{\eta\otimes\mu})(\al)= -\xi_{
    A^\ast}^{X_0,\,\nu\otimes\mu}(\al)-\lan\al,X_0\ran - 2\lan\al \wedge\d i_P\nu,\eta\ran,
    \quad\al\in\Omega^1(A).
\end{align*}

The previous equation is obviously equivalent to the duality
equation written in \cite{IglLopMarPad}. It  can also be rewritten
as
\begin{equation*}
P^\sharp(\xi_{A}^{\eta\otimes\mu})(\al)= -\xi_{
    A^\ast}^{\nu\otimes\mu}(\al)+(n-2)\lan\al,X_0\ran - 2\lan\al \wedge\d
    i_P\nu,\eta\ran,\quad \al\in\Omega^1(A),
\end{equation*}
or as
\begin{equation*}\label{relacao:quatro}
    P^\sharp(\xi_{A}^{{\phi_0},\,\eta\otimes\mu})(\al)= -\xi_{
    A^\ast}^{X_0,\,\nu\otimes\mu}(\al)+(n-2)\lan\al,X_0\ran - 2\lan\al \wedge\d
    i_P\nu,\eta\ran,\,  \al\in\Omega^1(A).
\end{equation*}

%%%%%%%%%%%%%%%%%%%%%%%  Modular class of Jacobi-Nijenhuis algebroids  %%%%%%%%%%%%%%%%%%%%%%%%%%%%%

\section{Modular classes of Jacobi-Nijenhuis algebroids}

Let $(A,{\phi_0})$ be a Jacobi algebroid and $N$ a Nijenhuis
operator. Consider a Jacobi bivector $P$ on $A$ compatible with the
Nijenhuis operator $N$. The sections $X_0=-P^\sharp({\phi_0})$ and
$X_1=-NP^\sharp({\phi_0})=-P^\sharp N^{\ast}({\phi_0})$ are
1-cocycles of the Lie algebroid $A^\ast_{N^\ast}$.

Since
 $(\hat A,\tilde P,N)$ is a Poisson-Nijenhuis Lie
 algebroid
it has a modular vector field (see \cite{raq}) given by
\begin{align*}
\hat X_{(N,\tilde P)}&=\xi_{\hat A^\ast_{N^\ast}}-N\xi_{\hat
A^\ast}\\
&=\hat\d_{\tilde P}(\tr N)=-\tilde P^\sharp(\hat \d\tr N)\\
&= -e^{-t}P^\sharp(d \tr N)=e^{-t}\d_P(\tr N).
\end{align*}

This $\hat A$-vector field is independent of the $Q_{\hat
A}$-section considered to compute the modular vector fields
$\xi_{\hat A^\ast_{N^\ast}}$ and $\xi_{\hat A^\ast}$. So  the
equation (\ref{relacao:modular:dualAum}) implies
\begin{align*}
   \hat X_{(N,\tilde
   P)}&=e^{-t}(\xi_{A^\ast_{N^\ast}}^{\nu\otimes\mu}-N\xi_{A^\ast}^{\nu\otimes\mu})
\end{align*}
and equation (\ref{relacao:modular:dualA}) implies
\begin{align*}
\hat X_{(N,\tilde
   P)}
   &=e^{-t}(\xi_{A^\ast_{N^\ast}}^{X_0,\,\nu\otimes\mu}-N\xi_{A^\ast}^{X_0,\,\nu\otimes\mu}),
\end{align*}
therefore
\begin{equation}%\label{}
\xi_{A^\ast_{N^\ast}}^{X_0,\,\nu\otimes\mu}-N\xi_{A^\ast}^{X_0,\,\nu\otimes\mu}=\d_P(\tr
N).
\end{equation}

 This relation motivates the next definition.

\begin{defi}
  The \textbf{modular vector field} of the  Jacobi-Ni\-je\-nhuis algebroid
  $(A,{\phi_0}, P, N)$ is defined by
    $$X_{(N,P)}=\xi_{A^\ast_{N^\ast}}-N\xi_{A^\ast}=\xi^{X_1}_{A^\ast_{N^\ast}}-N\xi^{X_0}_{A^\ast}$$
and is independent of the  section of $Q_A$ chosen. Its cohomology
class is called the \textbf{modular class} of $(A, {\phi_0}, P, N)$
and is denoted by ${\ds \m^{(N,P)} A=\brr{X_{(N,P)}}}$.
\end{defi}

\begin{rmk}
In fact, the modular class defined above is $\m({N^\ast})$, the
relative modular class of the Lie algebroid morphism $N^\ast:
A^\ast_{N^\ast}\ra A^\ast$ \cite{KosmannWeinstein}. As in the
Poisson case, $\m(N)$ and $\m({N^\ast})$  are related by $P$:
$$
P^\sharp\m({N})=-\m({N^\ast}).
$$
\end{rmk}

Following \cite{raq}, if $N$ is non-degenerated, we have a hierarchy
of  $\hat A$-vector fields:
$$
\hat X_{(N, \tilde P)}^{i+j}=N^{i+j-1}\hat X_{(N,\tilde
P)}=\d_{N^i\tilde P}h_j=\d_{N^j\tilde P}h_i,
$$
and a hierarchy of $A$-vector fields
$$
 X_{(N,  P)}^{i+j}=N^{i+j-1} X_{(N,
P)}=\d_{N^i P}h_j=\d_{N^jP}h_i,
$$
where
\begin{equation}
h_0=\ln (\det N) \quad \text{and} \quad
h_i=\frac{1}{i}\tr{N^i},\quad (i\neq 0, \quad i,j\in\Z).
\end{equation}

These hierarchies cover two  hierarchies, one on $M\times\R$ and
another one on $M$:

 The hierarchy on $M$ is given by
\begin{equation}
X^{i+j}_M=\rho
(X^{i+j}_{(N,P)})=-P_M^{i\,\sharp}(dh_j)=-P_M^{j\,\sharp}(dh_i)
\end{equation}
and the hierarchy on $M\times\R$ is given by
\begin{align*}
\hat X_{i+j}&=\hat\rho(\hat X_{(N, \tilde
P)}^{i+j})=\hat\rho(N^{i+j-1}\hat X_{(N,\tilde
P)})\\&=e^{-t}\hat\rho(X^{i+j}_{(N,P)})\\
&=e^{-t}\left(\rho (X^{i+j}_{(N,P)})+\lan {\phi_0}, X^{i+j}_{(N,P)}\ran \frac{\p}{\p t} \right)\\
&=e^{-t} \left(X_{M}^{i+j}+\lan \d h_j,
N^iP({\phi_0})\ran\frac{\p}{\p
t}\right)\\
&=e^{-t} \left(X_{M}^{i+j}+ \lan dh_j, E_M^i \ran\frac{\p}{\p t}\right)\\
&= e^{-t}\left(-(N^iP)_M^\sharp(dh_j)+\lan dh_j, E_M^i
\ran\frac{\p}{\p t}\right),
\end{align*}
where $((N^iP)_M,E_M^i)$ is the Jacobi structure  on $M$ induced by
the Jacobi  algebroid $(A,\phi_0,N^iP)$ (see
(\ref{jacobibivector:manifold:one}) and
(\ref{jacobibivector:manifold:two})).

This way we have proven the next theorem, which is a generalization
to Jacobi-Nijenhuis algebroids of the analogous result for
Poisson-Nijenhuis Lie algebroids \cite{raq} (see \cite{FerDam,
KosmannMagri0} for the  Poisson-Nijenhuis manifold case).

\begin{theorem}
\label{thm:main} Let $(A, {\phi_0}, P, N)$ be a Jacobi-Nijenhuis
algebroid with $N$ a non-degenerated Nijenhuis operator compatible
with $P$. Then the modular vector field $X_{(N, P)}$ is a
$\d_{NP}$-coboundary and determines a hierarchy of vector fields
\begin{equation}
\label{eq:A:hierarchy} X_{(N,
P)}^{i+j}=N^{i+j-1}X_{(N,P)}=\d_{N^iP}h_j=\d_{N^jP}h_i, \quad
(i,j\in\Z)
\end{equation}
where
\begin{equation}
\label{eq:Hamiltonians} h_0=\ln (\det N) \quad \text{and} \quad
h_i=\frac{1}{i}\tr{N^i},\quad (i\neq 0).
\end{equation}

This hierarchy covers a hierarchy of vector fields on $M$ given by
\begin{equation}
X_M^{i+j}=-(N^iP)_M^\sharp(dh_j)=-(N^jP)_M^\sharp(dh_i),
\end{equation}
and defines a hierarchy of vector fields on the Lie algebroid
$TM\times\R$ given by
\begin{equation}
Y^{i+j}=X_{M}^{i+j}+ \lan dh_j, E_M^i \ran\frac{\p}{\p t},
\end{equation}
 where $((N^iP)_M,E_M^i)$ are the Jacobi structures  on $M$
 induced by the
Jacobi bivectors $N^iP$ on $A$.
\end{theorem}

 \begin{rmk} Some remarks should be made at this point.
First, one should notice that even if $N$ is degenerated the
hierarchy exists but only for $i+j>1$, i.e.,
$$
X_{(N, P)}^{i+j}=\d_{N^iP}h_j=\d_{N^{j-1}P}h_{i+1}, \quad (0\leq
i<j, 1<j).
$$

In case $N$ is degenerated we can always consider a non-degenerated
Nijenhuis operator of the form $N+\lambda I$, $\lambda$ constant,
and we obtain the same algebra of commuting integrals.

It is also important to observe that although the hierarchy of
vector fields on $A$ is defined by a Nijenhuis operator, we may not
have a Nijenhuis operator on $M$ nor on $M\times\R$ that generates
neither one of the covered hierarchies.
 \end{rmk}

We will finish with a relation between  the sequence of modular
vector fields of the Jacobi-Nijenhuis algebroid and  the sequence of
modular vector fields of the Jacobi bialgebroid (in the sense of
\cite{IglLopMarPad}).

First recall the relation (\ref{relation:modular:com:Marrero}):
\begin{equation}\label{marreromodular}
    \mathcal{M}_{(A,{\phi_0},P)}^{\nu}(\al)=\xi_{A^\ast}^{\nu\otimes\mu}(\al)+X_0(\al)+i_P\d\al-\mathrm{div}_\mu(\rho\circ
    P^\sharp(\al)).
\end{equation}

Now we have
\begin{align*}
\mathcal{M}_{(A_N,\phi_1,P)}^{\nu}(\al)-&N\mathcal{M}_{(A,{\phi_0},P)}^{\nu}(\al)=
\xi_{A_N^\ast}^{\nu\otimes\mu}(\al)+X_1(\al)+i_{P}\d_N\al\\
&-\mathrm{div}_\mu(\rho_N\circ
    P^\sharp(\al))
    -N(\xi_{A^\ast}^{\nu\otimes\mu})(\al)-NX_0(\al)\\
    &-i_P\d N^\ast\al+\mathrm{div}_\mu(\rho\circ
    P^\sharp(N^\ast\al)))\\
    &= \lan\al,\d_P(\tr N)\ran+i_P \d_N\al-i_P\d N^\ast
    \al
\end{align*}
or equivalently, since $i_{NP}\d =i_P\d_N$,
\begin{align*}
\mathcal{M}_{(A_N,\phi_1,P)}^{\nu}(\al)-&N\mathcal{M}_{(A,{\phi_0},P)}^{\nu}(\al)=\\
    &=\lan\al,\d_P(\tr N)\ran+i_{NP} \d\al-i_P\d N^\ast\al\\
    &=\mathcal{M}_{(A,{\phi_0},NP)}^{\nu}(\al)-N\mathcal{M}_{(A,{\phi_0},P)}^{\nu}(\al).
\end{align*}

The vector field
$$
\mathcal{M}_{(N,P)}=\mathcal{M}_{(A_N,\phi_1,P)}^\nu-N\mathcal{M}_{(A,{\phi_0},P)}^\nu
$$
does not depend on the top-section of $A^\ast$ chosen and is related
with $X_{(N,P)}$ by
\begin{equation}%\label{}
  \lan
\al, \mathcal{M}_{(N,P)}\ran =\lan\al,X_{(N,P)}\ran+i_P
\d_N\al-i_P\d N^\ast \al.
\end{equation}

\

\begin{comment}
In particular, if $\al=\d f$, $f\in C^\infty(M)$, we have $i_P
\d_N\d f=-i_P\d N^\ast \d f$ and
$$
\rho(\mathcal{M}_{(N,P)})=-\rho(X_{(N,P)})+2i_P\d_N\df=-\rho(P^\sharp(\d\tr
N))=P_M(d\tr N).
$$
So $X_{N,P}$ and $\mathcal{M}_{(N,P)}$ cover the same hierarchy on
$M$. %Obviously, the  modular fields of the triangular Lie algebroids
%$(\hat A, \hat A^\ast,N^k\tilde P)$  also defines  the same
%hierarchy on $T(M\times\R)$ as $\hat X_{(N,\tilde P)}$.
\end{comment}

\begin{ex} Consider a Jacobi-Nijenhuis manifold $(M,  (\Lambda,E), \mathcal{N})$.
The modular class of the Jacobi manifold $(M, (\Lambda,E))$ is
 defined by
 (see \cite{Vaisman, IglLopMarPad})
$$
2\brr{V^{(\Lambda,E)}}=\!\!\!\mod(T^\ast M\times
\R)-(n+1)\brr{(E,0)}
$$
so
$$
\brr{V^{\mathcal{N}(\Lambda,E)}}-\mathcal{N}\brr{V^{(\Lambda,E)}}=\frac{1}{2}\brr{\d_*(\tr
\mathcal{N})}=\frac{1}{2}\!\!\!\!\mod^{(\mathcal{N},(\Lambda,E))}(T^\ast
M\times\R)
$$
and we have  the analogous relation as in the Poisson case.
\end{ex}

%%%%%%%%%%%%%%%%%%%%%  BIBLIOGRAPHY %%%%%%%%%%%%%%%%%%%%%%%

\end{document}